\documentclass{amsart}

\title{Desperately seeking mathematical proof}
\author{Melvyn B. Nathanson}
\address{Department of Mathematics,
Lehman College (CUNY),
Bronx, New York 10468, and 
CUNY Graduate Center, New York, NY 10016}
\email{melvyn.nathanson@lehman.cuny.edu}
\thanks{This paper elaborates one of the themes in Nathanson~\cite{nath08}.  I thank Saul Kripke for many useful discussions of various topics in the philosophy of mathematics.}

\begin{document}

\maketitle

How do we decide if a proof is a proof?  Why is a random sequence of true statements not a proof of the Riemann hypothesis?  That is, how do we know if a purported proof of a true theorem  is a proof?

Let's start with remarks on language.  The phrase ``true theorem'' is  redundant, and I won't use it again, since the definition of ``theorem'' is ``true mathematical statement,''  and the expression ``false theorem,'' like ``false truth,'' is contradictory.  I shall write ``theorem''  in quotation marks to denote a mathematical statement that is asserted to be true, that is, that someone claims is a theorem.   It is important to remember that a theorem is not false until it has a proof; it is only unproven until it has a proof.  Similarly, a false mathematical statement is not true until it has been shown that a counterexample exists or that it contradicts a theorem; it is just a ``statement'' until its untruth is demonstrated.  

The history of mathematics is full of philosophically and ethically troubling reports about bad proofs of theorems.  
For example, the fundamental theorem of algebra states that every polynomial of degree $n$ with complex coefficients has exactly $n$ complex roots.  D'Alembert published a proof in 1746, and the theorem became  known  ``D'Alembert's theorem,'' but the proof was wrong.  Gauss published his first proof of the fundamental theorem in 1799, but this, too, had gaps.  Gauss's subsequent proofs, in 1816 and 1849, were OK.   It seems to have been hard to determine if a proof of the fundamental theorem of algebra was correct.  Why?

Poincar\' e was awarded a prize from King Oscar II of Sweden and Norway for a paper on the three-body problem, and his paper was published in \emph{Acta Mathematica} in 1890.  But the published paper was not the prize-winning paper.  The paper that won the prize contained serious mistakes, and Poincar\' e and other mathematicians, most importantly, Mittag-Leffler, engaged in a conspiracy to suppress the truth and to replace the erroneous paper with an extensively altered and corrected one.

There are simple ways to show that a purported proof of a false mathematical statement is wrong.  For example, one might find a mistake in the proof, that is, a line in the proof that is false.  Or one might construct a counterexample to the  ``theorem.''   One might also be able to prove that the purported theorem is inconsistent with a known theorem.   I assume, of course, that mathematics is consistent.

To find a flaw in the proof of a theorem  is more complex, since no counterexample will exist, nor will the theorem contradict any other theorem.  A proof typically consists of a series of assertions, each leading more or less to the next, and concluding in the statement of the theorem.  How one gets from one assertion to the next can be complicated, since there are usually gaps.  We have to interpolate the missing arguments, or at least believe that a good graduate student or an expert in the field can perform the interpolation.  Often the gaps are explicit.  A typical formulation (Massey~\cite[p. 88]{mass91}) is:
\begin{quotation}
By the methods used in Chapter III, we can prove that the group $\pi(X)$ is characterized up to isomorphism by this theorem.  We leave the precise statement and proof of this fact to the reader. \end{quotation}
There is nothing improper about gaps in proofs, nor is there any reason to doubt that most gaps can be filled by a competent reader, exactly as the author intends.  The point is simply to emphasize that proofs have gaps and are, therefore,  inherently incomplete and sometimes wrong.  We frequently find statements like the following (Washington~\cite[p. 321]{wash97a}): 
\begin{quotation}
The Kronecker-Weber theorem asserts that every abelian extension of the rationals is contained in a cyclotomic field.  It was first stated by Kronecker in 1853, but his proof was incomplete \ldots.  The first proof was given by Weber in 1886 (there was still a gap \ldots).  
\end{quotation}

There is a lovely but probably apocryphal anecdote about Norbert Weiner.  Teaching a class at MIT, he wrote something on the blackboard and said it was ``obvious.''  One student had the temerity to ask for a proof.  Weiner started pacing back and forth, staring at what he had written on the board and  saying nothing.   Finally, he left the room, walked to his office, closed the door, and worked.  After a long absence he returned to the classroom.  ``It \emph{is} obvious,'' he told the class, and continued his lecture.

There is another reason why proofs are hard to verify:  Humans err.  We make mistakes and others do not necessarily notice our mistakes.  
Hume~\cite[Part IV, Section I]{hume1739}  expressed this beautifully in 1739:  
\begin{quotation}
There is no Algebraist nor Mathematician so expert in his science, as to place entire confidence in any truth immediately upon his discovery of it, or regard it as any thing, but a mere probability.  Every time he runs over his proofs, his confidence increases; but still more by the approbation of his friends; and is raised to its utmost perfection by the universal assent and applauses of the learned world.
\end{quotation}
This suggests an important reason why ``more elementary'' proofs are better than ``less elementary'' proofs:  The more elementary the proof, the easier it is to check and the more reliable is its verification.   We are less likely to err.  ``Elementary'' in this context does not mean elementary in the sense of elementary number theory, in which one tries to find proofs that do not use contour integrals and other tools of analytic function theory.  On elementary versus analytic proofs in number theory, I once wrote~\cite[p. ix]{nath00aa},
\begin{quotation}
In mathematics, when we want to prove a theorem, we may use any method.  The rule is ``no holds barred.''  It is OK to use complex variables, algebraic geometry, cohomology theory, and the kitchen sink to obtain a proof.  But once a theorem is proved, once we know that it is true, particularly if it is a simply stated and easily understood fact about the natural numbers, then we may want to find another proof, one that uses only ``elementary arguments'' from number theory.  Elementary proofs are not better than other proofs,\ldots.
\end{quotation}
I've changed my mind.  In this paper I argue that elementary (at least, in the sense of easy to check) proofs really are better.

Many mathematicians have the opposite opinion; they do not or cannot distinguish the beauty or importance of a theorem from its proof.   A theorem that is first published with a long and difficult proof is highly regarded.  Someone who, preferably many years later, finds a short proof is  ``brilliant.''  But if the short proof had been obtained in the beginning, the theorem might have been disparaged as an  ``easy result.''   Erd\H os was a genius at finding brilliantly simple proofs of deep results, but, until recently, much of his work was ignored by the mathematical establishment.

Erd\H os often talked about ``proofs from the Book.''   The ``Book'' would contain a perfect proof for every theorem, where a perfect proof was short, beautiful, insightful, and made the theorem instantaneously and obviously true.  We already know the ``Book proofs'' of many results.  I would argue that we do not, in fact,  fully understand a theorem until we have a proof that belongs in the Book.  It is impossible, of course, to know that  every theorem has a ``Book proof,'' but I certainly have the quasi-religious hope that all theorems do have such proofs. 

There are other reasons for the persistence of bad proofs of theorems.  Social pressure often hides mistakes in proofs.  In a seminar lecture, for example, when a mathematician is proving a theorem, it is technically possible to interrupt the speaker in order to ask for more explanation of the argument.  Sometimes the details will be forthcoming.  Other times the response will be that it's ``obvious'' or ``clear'' or ``follows easily from previous results.''  Occasionally speakers respond to a question from the audience with a look that conveys the message that the questioner is an idiot.  That's why most mathematicians sit quietly through seminars, understanding very little after the introductory remarks, and applauding politely at the end of a mostly wasted hour. 

Gel'fand's famous weekly seminar, in Moscow and at Rutgers, operated in marked contrast with the usual mathematics seminar or colloquium. One of the joys of Gel'fand's seminar was that he would constantly  interrupt the speaker to ask simple questions and give elementary examples.  Usually the speaker would not get to the end of his planned talk, but the audience would actually learn some mathematics.

The philosophical underpinning to this discussion is the belief that  ``mathematical'' objects exist in the real world, and that mathematics is the science of discovering and describing their properties, just as ``physical'' objects exist in the real world and physics is the science of discovering and describing their properties.  This is in contrast to an occasionally  fashionable notion that mathematics is the logical game of deducing conclusions from interesting but arbitrarily chosen finite sets of axioms and rules of inference.  If the mathematical world is real, then it is unlikely that it can be encapsulated in any finite system.  There are, of course, masterpieces of mathematical exposition that develop a deep subject from a small set of axioms.  Two examples of  books that are perfect in this sense are Weil's \emph{Number Theory for Beginners}, which is, unfortunately, out of print, and Artin's \emph{Galois Theory}. 
Mathematics can be done scrupulously.  

Different theorems can be proven from different assumptions.  The compendium of mathematical knowledge, that is, the collection of theorems, becomes a social system with various substructures, analogous to clans and kinship systems, and a newly discovered theorem has to find its place in this social network.  To the extent that a new discovery fits into an established community of mathematical truths, we believe it and tend to accept its proof.  A theorem that is an ``outsider,'' a kind of  social outlaw, requires more rigorous proof, and finds acceptance harder.  

Wittgenstein~\cite[p. 401]{witt83} wrote, 
\begin{quotation}
If a contradiction were now actually found in arithmetic -- that would only prove that an arithmetic with \emph{such} a contradiction in it could render very good service; and it would be better for us to modify our concept of the certainty required, than to say it would really not yet have been a proper arithmetic.
\end{quotation}
This passage (still controversial in the philosophy of mathematics) evidences a pragmatic approach to mathematics that describes how mathematicians behave in the privacy of their offices, in contrast to our more pietistic public pronouncements.

Perhaps we should discard the myth that mathematics is a rigorously deductive enterprise.  It may be more deductive than other sciences, but hand-waving is intrinsic.  We try to minimize it and we can sometimes escape it, but not always, if we want to discover new theorems.

\def\cprime{$'$} \def\cprime{$'$} \def\cprime{$'$}
\providecommand{\bysame}{\leavevmode\hbox to3em{\hrulefill}\thinspace}
\providecommand{\MR}{\relax\ifhmode\unskip\space\fi MR }
\providecommand{\MRhref}[2]{%
  \href{http://www.ams.org/mathscinet-getitem?mr=#1}{#2}
}
\providecommand{\href}[2]{#2}

\end{document}